\documentclass{article}
\usepackage{mathrsfs}
\usepackage{amsthm}
\usepackage{fullpage}
\usepackage{enumerate}
\usepackage{amsmath}     
\usepackage{amssymb}     
\usepackage{amsfonts} 
\usepackage{hyperref}

\newtheorem{theorem}{Theorem}
\newtheorem{lemma}{Lemma}
\newtheorem{proposition}{Proposition}
\theoremstyle{definition}
\newtheorem{definition}{Definition}[section]
\newtheorem{corollary}{Corollary}

\begin{document}

\title{Optimal domains of Volterra-type operators on Bergman spaces}
\author{Zengjian Lou and Jinbang Wen}
\date{}

\maketitle
\begin{abstract}
\large We consider bounded Volterra-type operators $T_g$ from the weighted Bergman space $A^q_{\beta}$ to $A^p_{\alpha}$ for non-constant $g\in\operatorname{Hol}(\mathbb D)$.\ The holomorphic optimal domain $[T_{g},A^p_{\alpha}]$ is introduced and shown to be a quasi-Banach space.\ We show that the intersection of the optimal domains is equal to $A^q_{\beta}$ if $p\le q$, while if $p>q$, it is strictly larger than $A^q_{\beta}$.\ Moreover, $A^q_{\beta} \subsetneq[T_{g},A^p_{\alpha}]$ is established for $p\neq q$.\ If $p=q$, by contrast, $A^p_{\beta}$ may coincide with $ [T_{g},A^p_{\alpha}]$.\ The multiplier algebras and strict inclusions with respect to varying parameters are established.\ For example, we determine exactly when $[T_{g},A^s_{\gamma}] \subset [T_{g},A^p_{\alpha}]$ for some $g\in \operatorname{Hol}(\mathbb{D})$.\ We further characterize the coincidence of the holomorphic optimal domains.


\end{abstract}
\section{Introduction}
\large
\quad Let $\mathbb{D}$ denote the unit disk of the complex plane $\mathbb{C}$, $\mathbb{T}$ denote its boundary.\ We denote as $\operatorname{Hol}(\mathbb{D})$ the spaces of holomorphic functions on $\mathbb{D}$.

For a fixed function $g\in\operatorname{Hol}(\mathbb{D})$, we consider the Volterra-type operator acting on holomorphic functions,
\[
T_g(f)(z)=\int_{0}^{z} f(\zeta) g^{\prime}(\zeta)\,d\zeta,\qquad z\in\mathbb{D},\quad f\in\operatorname{Hol}(\mathbb{D}).
\]
The operator $T_g$ specializes to the Volterra operator when $g(z)=z$ and to the Cesàro operator when $g(z)=-\log(1-z)$.\ We also consider the integral operator
\[
S_{g}(f)(z)=\int_{0}^{z}f'(\zeta)g(\zeta)\, d\zeta,\quad z\in \mathbb{D},
\]
and the multiplication operator $M_{g}(f)(z)=f(z)g(z)$.\ We denote by $\mathcal M(X,Y)$ the space of multipliers $g$ such
that $M_g:X\to Y$ is bounded.\ The principal objective in the investigation of $T_g$ is its action on spaces of holomorphic functions, with particular emphasis on the Hardy spaces and Bergman spaces on $\mathbb{D}$.\ $H^\infty$ denotes the space of bounded holomorphic functions on $\mathbb{D}$, endowed with the supremum norm.\ We also consider the weighted Bergman space $A^p_\alpha$ on $\mathbb D$ for $0< p<\infty$ and $-1<\alpha<\infty$
$$
A^p_\alpha=\left\{f\in \operatorname{Hol}(\mathbb{D}):\,\|f\|^p_{p,\alpha} =(\alpha+1)\int_\mathbb{D} |f(z)|^p(1-|z|^2)^\alpha \, dA(z) < \infty\right\},$$
where $dA(z)$ denotes the normalized area measure on $\mathbb{D}$, and we set $dA_\alpha(z)=(\alpha+1)(1-|z|^2)^\alpha dA(z)$.\ When $\alpha=0$, this yields the classical Bergman spaces $A^p$.\ See \cite{13} for more information on Bergman spaces.

Pommerenke \cite{18} established that the set of symbols $g$ for which $T_g$ is bounded on the Hilbert space $H^2$ is precisely $BMOA$.\ Subsequently, Aleman and Siskakis \cite{3} extended this result to weighted Bergman spaces, proving that, for all $1\leq p<\infty$ and $-1<\alpha<\infty$, $T_g$ is bounded on $A^p_{\alpha}$ if and only if $g\in\mathcal{B}$.\ Finally, Xiao \cite{21} characterized the boundedness of $T_g\colon A^q_\beta\to A^p_\alpha$ for all $0<p,q<\infty$ and $-1<\alpha,\beta<\infty$ as follows:\ When $q > p$, $g$ satisfies  
\[
\int_{\mathbb{D}} \Bigl( |\,g'(z)| \,(1-|z|^2)^{1-\frac{\beta}{q}+\frac{\alpha}{p}} \Bigr)^{\frac{pq}{q-p}} \, dA(z) < \infty;
\]  
when $q \le p$, $g$ satisfies  
\[
\sup_{z\in\mathbb{D}} |\,g'(z)| \,(1-|z|^2)^{1-\frac{2+\beta}{q}+\frac{2+\alpha}{p}} < \infty.
\]
Firstly, for $\tau\ge 0$, $\mathcal{B}^{\tau}$ is defined as
\[
\mathcal{B}^{\tau}= \left\{ g \in \operatorname{Hol}(\mathbb{D}) : \|g\|_{\mathcal{B}^{\tau}} = |g(0)| + \sup_{z \in \mathbb{D}} (1 - |z|^2)^{\tau} |g'(z)| < \infty\right\}.
 \]
With this norm $\mathcal B^\tau$ is a Banach space.\ For $\tau=1$ it coincides with the classical Bloch space $\mathcal{B}$.\ Moreover, the little Bloch space of $\mathbb D$, denoted by $\mathcal B_0$, is the closed subspace of $\mathcal B$ consisting of those functions $f$
for which
\[
\lim_{|z|\to 1^-}(1-|z|^2)f'(z)=0.
\]
 
Inspired by \cite{20}, we consider $V_{q,p} := \{g \in \operatorname{Hol}(\mathbb{D}) : T_g \text{ is bounded from } A_\beta^q \text{ to } A_\alpha^p\}$.\ These spaces are completely determined by the values of $p, \alpha$ and $q, \beta$:
\begin{itemize}
\item $q > p$: $V_{q,p} = W^r_\delta$;
\item $q \le p$: $V_{q,p} = \mathcal{B}^{\tau}\,(\tau=1-\frac{2+\beta}{q}+\frac{2+\alpha}{p})$.
\end{itemize}
Here $W^r_\delta:= \{g \in \operatorname{Hol}(\mathbb{D}) : g'\in A^r_\delta, \delta=\frac{pq}{q-p}(1-\frac{\beta}{q}+\frac{\alpha}{p}), r= \frac{pq}{q-p}\}$ with the quasi-norm
\[
\|g\|_{W^r_\delta} = |g(0)| + \int_\mathbb{D} |g'(z)|^r(1-|z|^2)^\delta \, dA(z)  < \infty
\]
is a quasi-Banach space.\ In particular, the parameters in this paper are chosen such that $\frac{2+\beta}{q}-\frac{2+\alpha}{p}\le 1$ for $q\le p$ or $\frac{1+\beta}{q}-\frac{1+\alpha}{p}< 1$ for $q> p$, otherwise $g$ is constant.

 For $0\le \kappa<\infty$, we also consider the Korenblum growth spaces \cite{15}, defined as
 \[
     \mathcal{A}^{-\kappa}=\left\{ f \in \operatorname{Hol}(\mathbb{D}) : \|f\|_{\mathcal{A}^{-\kappa}} = \sup_{z \in \mathbb{D}} (1 - |z|^2)^\kappa |f(z)| < \infty\right\}.
 \]For $\kappa=0$, we recover the classical $H^\infty$ space.

In \cite{8} and \cite{9}, Curbera and Ricker introduced the holomorphic optimal domain of the classical Ces\`aro operator on the Hardy spaces $H^p$.\ Subsequently, Bellavita et al.\ \cite{7} studied the optimal domain of the Volterra-type operator $T_g$ acting on $H^p$ and gave a full description of its structural properties.\ In the recent article \cite{2} the authors extended these results to the optimal domain problem for the operator $T_g$ on the Korenblum spaces.\ For other spaces, see \cite{11} and \cite{17}.
\begin{definition}
Let $0<p,q<\infty$, $-1<\alpha,\beta<\infty$, and $g\in V_{q,p}$.\ The optimal domain of $T_g$ on $A^p_{\alpha}$ is defined as
\[
[T_g, A^p_{\alpha}] = \bigl\{ f \in \operatorname{Hol}(\mathbb{D}) : T_g(f) \in A^p_{\alpha} \bigr\}
\]
with the natural quasi-norm $\|f\|_{[T_g, A^p_{\alpha}]} = \|T_g f\|_{p,\alpha}$.
\end{definition}
It is clear that $ [T_{g},A^p_{\alpha}] $ is a linear subspace of $ \operatorname{Hol}(\mathbb{D})$.\ If $g$ is constant, then trivially $ [T_{g},A^p_{\alpha}]= \operatorname{Hol}(\mathbb{D})$. 

In what follows, we study some basic properties of this quasi-normed space.\ However, unlike the Hardy space case, by reference \cite{4}, $T_g$ may be bounded below.\ To study the case where  $A^p_{\beta}=[T_{g},A^p_{\alpha}]$, we need the following definition \cite{24}.

\begin{definition}
For $z,w\in\mathbb D$, let
$\varphi_z(w)=\frac{z-w}{1-\overline z w}$.\
This induces the pseudohyperbolic metric
$\rho(z,w)=|\varphi_z(w)|$.\ For $0<r<1$, the pseudohyperbolic disk with center $z$ and radius $r$ is defined by
\[
\Delta(z,r)=\{w\in\mathbb D:\rho(z,w)<r\}.
\]
\end{definition}

Whereas in \cite{4} the lower estimate for $T_g$ on weighted Bergman spaces relies on Luecking’s inequalities, which are valid only for functions in the Bergman space, our approach yields the corresponding lower estimate for all holomorphic functions $f$ on $\mathbb{D}$~(see Proposition~\ref{pro}).\ As a consequence, under the standing assumptions we obtain the identification $A^p_\beta=[T_g,A^p_\alpha]$.

A natural and interesting problem is to characterize those functions $h\in V_{q,p}$ for which $[T_h,A^p_{\alpha}]=[T_g,A^p_{\alpha}]$.\ To address this question, we introduce the following spaces.

\begin{definition}
Let $0<p,q<\infty$, $-1<\alpha,\beta<\infty$, and $g\in V_{q,p}$.\ Define
\begin{align*}
V_g^p &=\left\{h\in V_{q,p}:[T_h,A^p_{\alpha}]\supseteq[T_g,A^p_{\alpha}]\right\}.
\end{align*}
\end{definition}

Let $0<p,q,s<\infty$ and $-1<\alpha,\beta,\gamma<\infty$.\ This paper is devoted to spaces $[T_g,A^{p}_{\alpha}]$, and $V_g^p$, with emphasis on their structural features.\ Our analysis sharpens the description of the holomorphic optimal domain $[T_g,A^{p}_{\alpha}]$ and the embedding problem.\ Beyond the Hardy-space setting, we establish for $g\in V_{q,p}$ several results on the optimal domain of $T_g\colon A^{q}_{\beta}\to A^{p}_{\alpha}$.\ We also study inclusion relations between $[T_g,A^{p}_{\alpha}]$ and $[T_g,A^{s}_{\gamma}]$ for some $g\in \operatorname{Hol}(\mathbb{D})$.\ We further determine, for $g_1,g_2\in V_{q,p}$, when the optimal domains $[T_{g_1},A^{p}_{\alpha}]$ and $[T_{g_2},A^{p}_{\alpha}]$ coincide.

Although certain theorems concerning the optimal domain in the Hardy space can be transferred to the weighted Bergman spaces by standard methods, whereas for Theorem~\ref{theo:2}, the generalization is achieved by a special property of growth spaces when $q\le p$, and by the multiplier properties between weighted Bergman spaces when $q>p$.\ The proof of Theorem~\ref{theo:5} makes use of the test functions, pseudo-hyperbolic disks and an $r$-lattice.\ The multiplier space properties play an essential role in the proof of Theorem~\ref{theo:7}.\ By means of an embedding of the same type as that for $A^p_{\alpha}$, we obtain Theorem~\ref{theo:8}.\ For the space $V_g^p$, we adopt new methods to relax the conditions required in the corresponding Hardy-space theorems, obtaining improved results---namely Theorem~\ref{theo:9}--\ref{theo:11} and Proposition~\ref{pro:4}.\ Finally, when $g\in \mathcal{B}$, by exploiting multiplier spaces properties, we provide several examples of $A^p_{\alpha}\subsetneq[T_{g},A^p_{\alpha}]$, which ultimately lead to a necessary condition for $T_g$ to be bounded below on $A^p_{\alpha}$.

Throughout the paper, for real-valued functions $f$ and $g$, we write $f\lesssim g$ whenever $f\le Cg$ for some constant $C>0$ that does not depend on the variables of $f$ and $g$;\ the constant $C$ is allowed to change from one occurrence to another.\ The relation $f\gtrsim g$ is understood in the analogous sense, and $f\sim g$ means that $f\lesssim g$ and $f\gtrsim g$ hold at the same time.\ For $h\in\operatorname{Hol}(\mathbb D)$, the symbol $\mathcal Z(h)$ denotes the zero set of $h$, with zeros counted according to multiplicity.\ Moreover, we write $p \in \mathcal{X} \setminus \mathcal{Y}$ to express that $p$ belongs to $\mathcal{X}$ and fails to belong to $\mathcal{Y}$.

\section{Main results and plan of the paper}

\quad From now on, unless otherwise specified, we always assume that $0<p,q<\infty$, $-1 < \alpha, \beta <\infty$, $g \in V_{q,p}$
is non-constant, and $\frac{2+\beta}{q}-\frac{2+\alpha}{p}\le 1$ for $q\le p$ or $\frac{1+\beta}{q}-\frac{1+\alpha}{p}< 1$ for $q> p$.

\begin{theorem}
\label{theo:1}
Equipped with the quasi-norm $\|\cdot\|_{[T_g,A^p_\alpha]}$, the space
$[T_g,A^p_\alpha]$ is a quasi-Banach space.\ It is further a Banach space whenever $p\ge1$.
\end{theorem}

 We note that the space $[T_{g}, A^p_{\alpha}]$ contains $A^q_{\beta}$.\ Along the same lines, the following theorem is proved:

\begin{theorem}
\label{theo:2}
\begin{enumerate}
\item[\textup{(a)}] If $p \leq q$, then
\[
A^q_{\beta}=\bigcap_{g\in V_{q,p}}\left[T_{g},A^p_{\alpha}\right].
\]
\item[\textup{(b)}]  If $ q < p$, then 
\[
A^q_{\beta} \subsetneq \bigcap_{g\in V_{q,p}}\left[T_{g},A^p_{\alpha}\right].
\]
\end{enumerate}
\end{theorem}
\begin{theorem}
\label{theo:3}
Let $ p \neq q$.\ Then $ A^q_{\beta} $ is strictly contained in $ [T_{g},A^p_{\alpha}]$.
\end{theorem}
However, the case $p=q$ is more delicate.\ In this case, $A^p_{\beta}$ and $ [T_{g},A^p_{\alpha}]$ are allowed to be equal.
\begin{theorem}
\label{theo:5}
Let $p=q$.\ If $T_g$ is bounded below, then $A^p_{\beta} =[T_{g},A^p_{\alpha}]$;\ otherwise, $A^p_{\beta} \subsetneq[T_{g},A^p_{\alpha}]$.
\end{theorem}

Obviously, the optimal domain is also a weighted Bergman space.\ Hence all the complications inherent to those spaces carry over to the analysis of $[T_g, A^p_\alpha]$.\ The question of characterizing those symbols $g$ for which the polynomials are dense in $[T_g, A^p_\alpha]$ remains open.\ We present explicit examples demonstrating  polynomials are not dense in $[T_{g},A^p_{\alpha}]$.\ Similar to the weighted Bergman spaces, we prove that the multiplier algebra remains unchanged.

\begin{theorem}
\label{theo:6}
The multiplication operator
\[
M_{h}:\left[T_{g},A^p_{\alpha}\right]\to\left[T_{g},A^p_{\alpha}\right]
\]
is bounded if and only if $h\in H^{\infty}$.\ 
\end{theorem}

It is well known \cite[Theorem~C]{5} that for weighted Bergman spaces, varying parameters lead to distinct inclusion relations.\ We establish that analogous inclusion relations also hold for the associated optimal domains of weighted Bergman spaces.
\begin{theorem}
\label{theo:7}
Let $0<p_1 < p_2 < \infty$.\ Assume additionally that $g\in V_{q,p_2}$ is non-constant satisfies: $T_g$ is bounded from $A_\beta^{q}$ to $A_\alpha^{p_2}$.\ Then we have that
 \begin{enumerate}
    \item[\textup{(a)}] If $q \le p_1$, then\[
[T_{g}, A^{p_2}_{\alpha}] \subsetneq [T_{g}, A^{p_1}_{\alpha}].
\]
\item[\textup{(b)}] If $q>p_1$ and $ \frac{2+\alpha}{p_2} - \frac{1}{q}< \frac{1+\alpha}{p_1}$, then
\[
[T_{g}, A^{p_2}_{\alpha}] \subsetneq [T_{g}, A^{p_1}_{\alpha}].
\]
\item[\textup{(c)}] The inclusion is compact.
\end{enumerate}
\end{theorem}
In Theorem~\ref{theo:8}, the definition of the optimal domain may be extended to allow $g\in \operatorname{Hol}(\mathbb{D})$ is non-constant;\ the corresponding optimal domain remains a quasi-Banach space.
\begin{theorem}
\label{theo:8}
Let $0<s<\infty$, $-1<\gamma<\infty$, and let
$g\in\operatorname{Hol}(\mathbb{D})$ be non-constant and such that $g'$ has only finitely many zeros in $\mathbb{D}$.\ Then 
\begin{enumerate}
\item[\textup{(a)}] If $s \leq p$, then $[T_{g},A^s_{\gamma}] \subset [T_{g},A^p_{\alpha}]$ if and only if
    \begin{equation*}
        \frac{2+\gamma}{s} \leq \frac{2+\alpha}{p}. 
    \end{equation*}
    
\item[\textup{(b)}] If $s > p$, then $[T_{g},A^s_{\gamma}] \subset [T_{g},A^p_{\alpha}]$ if and only if
    \begin{equation*}
        \frac{1+\gamma}{s} < \frac{1+\alpha}{p}. 
    \end{equation*}
    
\item[\textup{(c)}] Embeddings in the case $s > p$ are always compact, while embeddings for $s \leq p$ are compact if and only if strict inequality holds.
\end{enumerate}
\end{theorem}

We now consider $V^p_g$.\ The question is to characterize when two symbols $g_1,g_2\in V_{q,p}$ give rise to the same optimal domain, i.e., $[T_{g_1},A^p_{\alpha}] = [T_{g_2},A^p_{\alpha}]$.

For a finite or infinite sequence $\{z_n\}$ in $\mathbb D$ satisfying the Blaschke condition
$\sum_n(1-|z_n|)<\infty$, an associated Blaschke product is given by \[B(z)=\lambda z^m\prod_{n:\,z_n\ne0}\frac{|z_n|}{z_n}\frac{z_n-z}{1-\overline{z_n}z},
\]
where $|\lambda|=1$ and $m$ is the multiplicity of $0$ in $\{z_n\}$.\ A sequence $\{z_n\}$ in $\mathbb D$ is called uniformly separated if
\[
\inf_n\prod_{k\ne n}\rho(z_k,z_n)>0.
\]
\begin{theorem}
\label{theo:9}
Suppose that $g'=Bg'_1$, where $g'_1\in\operatorname{Hol}(\mathbb{D})$ is zero-free and $B$ is a Blaschke product whose zero sequence, counted with multiplicity, is a finite union of uniformly separated sequences.\ Then
\[
V_g^p=T_{g_1}(H^\infty)+\mathbb{C}.
\]
\end{theorem}

\begin{theorem}
\label{theo:11}
Let $g'$ be a holomorphic function on $\bar{\mathbb{D}}$.\ Suppose that $\lambda_j\in\overline{\mathbb{D}}\ (j=1,2,...,N)$ are the points such that $g^{\prime}(\lambda_j)=0$.\ Then the following hold:
\begin{enumerate}
\item[\textup{(a)}] If $|\lambda_j|\neq 1$ for every $j$, $[T_g,A^p_{\alpha}]=[T_z,A^p_{\alpha}]=A^p_{\alpha+p}$.
\item[\textup{(b)}] If there exists an index $j$ such that $|\lambda_j|=1$, then $A^p_{\alpha+p}\subsetneq[T_g,A^p_{\alpha}]$.
\end{enumerate}
\end{theorem}This theorem shows that there exists a symbol $g$ for which
$[T_g,A^p_\alpha]$ strictly contains $A^p_{\alpha+p}$.

The paper proceeds as follows.\ Section~3 is devoted to the structure of the optimal domain and to the inclusion $A^q_{\beta} \subset[T_{g},A^p_{\alpha}]$.\ It contains the proofs of Theorems~\ref{theo:1}-\ref{theo:5}.\ Section~4 addresses the pointwise multipliers of the optimal domain and the embedding problem between different optimal domains.\ Moreover, paralleling the situation for weighted Bergman spaces, we establish inclusion relations among optimal domains with respect to the parameter and prove Theorems~\ref{theo:7} and \ref{theo:8}.\ In Section 5, the structure of $V^p_g$ is sketched, and a characterization is given of when two holomorphic optimal domains are identical for special choices of $g$.\ In the last section, several examples of the proper inclusion $A^p_\alpha\subsetneq[T_g,A^p_\alpha]$ are given with $g\in \mathcal{B}$, which leads to Corollary~\ref{cor:3}.

\section{The structure of the holomorphic optimal domain}

\quad $[T_g,A^p_\alpha]$ is a linear space follows directly from its definition.\ Our main task is therefore to confirm that it forms a quasi‑Banach space.\ As noted in \cite{4}, the injectivity of $T_g$ for any non‑constant $g\in\operatorname{Hol}(\mathbb{D})$ guarantees that $\|\cdot\|_{[T_g,A^p_\alpha]}$ is a genuine quasi‑norm for all $p>0$.

In the following proofs we will repeatedly use the equivalent description of the $A^p_{\alpha}$ norm \cite[Theorem~4.28]{24}: for $p>0$, and any $f\in\operatorname{Hol}(\mathbb{D})$ with $f(0)=0$, one has $f\in A^p_{\alpha}$ if and only if $f'\in A^p_{\alpha+p}$, and $\|f\|_{p,\alpha}\sim\|f'\|_{p,\alpha+p}$ where the implied constants depend only on $p$ and~$\alpha$.

\begin{proof}[\textbf{Proof of Theorem~\ref{theo:1}}]
The proof is similar to that of Theorem~1 in \cite{7} for Hardy spaces.\ We begin by showing that, on every compact subset $K$ of the unit disk, the point evaluation
functionals of $[T_g,A^p_\alpha]$ are uniformly bounded.

Fix $z_0\in\mathbb D$, for $h\in A^p_{\alpha}$, from \cite[Proposition 1.1]{13} we infer that $h'(z)$ belongs to $A^p_{\alpha+p}$, and
one can find a constant $C>0$ such that
\[
|h'(z_0)|\le\frac{C}{(1-|z_0|)^{1+\frac{\alpha+2}{p}}}\|h\|_{p,\alpha}.
\]
Let $f\in[T_g,A^p_\alpha]$.\ Then
\[
|f(z_0)g'(z_0)|=|T_g(f)'(z_0)|
\le\frac{C}{(1-|z_0|)^{1+\frac{\alpha+2}{p}}}\|f\|_{[T_g,A^p_\alpha]}.
\]
Consequently, if $g'(z_0)\neq0$, we can divide by $|g'(z_0)|$, and the continuity of $g'$
guarantees that the point evaluations are bounded in a neighbourhood of $z_0$.\ If $g'(z_0)=0$, since the zeros of $g'$ are isolated, there is an annulus $N_{z_0}$ centred at $z_0$ such that
$$
\min_{\zeta\in N_{z_0}}|g'(\zeta)|>\varepsilon>0 .
$$
The subharmonicity of $|f|$ and the sub-mean value property give
$$
\begin{aligned}
|f(z_0)|
&\le\frac{1}{|N_{z_0}|}\int_{N_{z_0}}|f(w)|\,dA(w)\\
&\le\frac{1}{|N_{z_0}|}\int_{N_{z_0}}
\frac{C}{(1-|w|)^{1+\frac1p}\displaystyle\min_{\zeta\in N_{z_0}}|g'(\zeta)|}\,
\|f\|_{[T_g,A^p_\alpha]}\,dA(w)\\
&=C_1\|f\|_{[T_g,A^p_\alpha]},
\end{aligned}
$$
where $C_1>0$ depends only on $z_0$.\ By continuity of $f$, the point evaluations are uniformly bounded also in a neighbourhood of $z_0$.\ Because $K$ is compact, finitely many such neighbourhoods cover $K$, and hence the point evaluations are uniformly bounded on $K$.

We conclude by proving the completeness.\ Let $\{f_n\}_{n\ge1}\subset[T_g,A^p_\alpha]$ be a Cauchy sequence.\ Then $T_g(f_n)$ is a Cauchy sequence in $A^p_\alpha$, which is complete, so as $n\to \infty$, $T_g(f_n)\to u$ in $A^p_\alpha$ for some $u$.\ By the uniform boundedness of the point evaluations above, $f_n\to f$ locally uniformly for some $f\in \operatorname{Hol}(\mathbb{D})$, and passing to limits gives $u'=g'f$, so $u=T_g(f)$, which completes the proof.\ Since $T_g$ is
injective for non-constant $g$, $[T_g,A^p_\alpha]$ is a Banach space for $p\ge1$, and a quasi-Banach space for $0<p<1$.
\end{proof}
To prove Theorem~\ref{theo:2}, we first record the following lemma see \cite[Theorem~C]{5}.\ It
will be used repeatedly in the sequel.
\begin{lemma}\label{0} 
Suppose $0 < p, q < \infty$ and $-1<\alpha, \beta <\infty$.\ Then
\begin{enumerate}
\item[\textup{(a)}] If $q \leq p$, then $A_{\beta}^{q} \subset A_{\alpha}^{p}$ if and only if
    \begin{equation*}
       \frac{2+\beta}{q}\leq \frac{2+\alpha}{p}.
    \end{equation*}
    
\item[\textup{(b)}] If $q > p$, then $A_{\beta}^{q} \subset A_{\alpha}^{p}$ if and only if
    \begin{equation*}
       \frac{1+\beta}{q}< \frac{1+\alpha}{p}.
    \end{equation*}
    
\item[\textup{(c)}] Embeddings in the case $q > p$ are always compact, while embeddings for $q \leq p$ are compact if and only if strict inequality holds.
\end{enumerate}
\end{lemma}
\begin{proof}[\textbf{Proof of Theorem~\ref{theo:2}}]
Since $ T_{g}(A^q_{\beta})\subset A^p_{\alpha}$ for every $ g\in V_{q,p} $, we note that
$$ 
A^q_{\beta}\subseteq\bigcap_{g\in V_{q,p}}\left[T_{g},A^p_{\alpha}\right].
$$ 

We first prove that: $f\in \bigcap_{g\in V_{q,p}} [T_g, A^p_{\alpha}]$ if and only if $S_f$ is a bounded operator from $V_{q,p}$ to $A^p_{\alpha}$.\ Let $f\in \bigcap_{g\in V_{q,p}} [T_g, A^p_{\alpha}]$.\ The identity $T_g(f) = S_f(g)$ together with the choice of $f$ shows that $S_f(g) \in A^p_{\alpha}$ for every $g\in V_{q,p}$.\ An application of the closed graph theorem (for quasi-Banach spaces if necessary;\ see \cite{19}) then yields the boundedness of $S_f : V_{q,p} \to A^p_{\alpha}$.\ If $S_f$ is bounded, then $T_g(f)=S_f(g)\in A^p_{\alpha}$ for all $g\in V_{q,p}$, the converse is true.

For $q > p$, we show that the boundedness of $S_f : V_{q,p} \to A^p_{\alpha}$ is equivalent to $ f\in A^q_{\beta}$.\ If $S_f$ is bounded, then it follows from the definitions of $S_f$ and $V_{q,p}$ that $f$ is a multiplier from $A_\delta^r$ to $A_{\alpha+p}^{p}$, where $
r=\frac{qp}{q-p}>p$ and $\delta=(1-\frac{p}{q}+\frac{\alpha}{p})\frac{qp}{q-p}$.\ By \cite[Theorem~1]{23}, this multiplier condition is equivalent to $f\in A_\gamma^s$, where $s=\frac{rp}{r-p}$ and $\frac{\gamma}{s}=\frac{\alpha+p}{p}-\frac{\delta}{r}$.\ Substituting the expressions for $r$ and $s$ yields
$$s=\frac{\frac{pq}{q-p}\cdot p}{\frac{pq}{q-p}-p}=q, \qquad\text{and then}\qquad \gamma=q(1+\frac{\alpha}{p}-(1-\frac{\beta}{q}+\frac{\alpha}{p}))=\beta.$$
That is $f\in A_\beta^q$.\ On the other hand, by the closed graph theorem, $S_f$ is bounded for $f \in A^q_\beta$.
 
In the case $q \leq p$, we show that the boundedness of $S_f : V_{q,p} \to A^p_{\alpha}$ is equivalent to $ f\in A^p_{\frac{p(2+\beta)}{q}-2}$.\ To prove this, set $\tau=1-\frac{2+\beta}{q}+\frac{2+\alpha}{p}$.\ By the assumption of the theorem, $\tau\ge 0$.\ Let $f \in A^p_{\frac{p(2+\beta)}{q}-2}$.\ For every $g \in V_{q,p}$, we have 
\[
\| S_{f}(g) \|_{p,\alpha} \sim \| g' f \|_{p,\alpha+p} \lesssim \| f \|_{p,\frac{p(2+\beta)}{q}-2} \| g \|_{\mathcal{B}^\tau}.
\]
That is $S_{f}$ is bounded.

Conversely, if $S_{f}: V_{q,p} \to A^p_{\alpha}$ is bounded.\ From \cite[Theorem~1.2]{1} and \cite[Lemma 3.1]{22}, we obtain that there exist two holomorphic functions $g_{1},g_{2}\in \mathcal{B}^\tau$ such that
\[
|g_{1}'(z)| + |g_{2}'(z)| \geq \frac{1}{(1-|z|^{2})^\tau} \quad \text{for all } z \in \mathbb{D}.
\]
Therefore, $S_f(g_1) \in A^p_{\alpha}$, $S_f(g_2) \in A^p_{\alpha}$, and
\begin{align*}
\infty &> \int_{\mathbb{D}} |S_f(g_1)|^p (1-|z|^2)^{\alpha}\, dA(z) 
          + \int_{\mathbb{D}} |S_f(g_2)|^p (1-|z|^2)^{\alpha}\, dA(z) \\
&\sim \int_{\mathbb{D}} \bigl( |g_{1}^{\prime}(z)|^{p} + |g_{2}^{\prime}(z)|^{p} \bigr) |f(z)|^{p} (1-|z|^2)^{\alpha+p} \, dA(z) \\
&\geq \frac{1}{2^{p}} \int_{\mathbb{D}} \bigl( |g_{1}^{\prime}(z)| + |g_{2}^{\prime}(z)| \bigr)^{p} |f(z)|^{p} (1-|z|^2)^{\alpha+p} \, dA(z) \\
&\geq \frac{1}{2^{p}} \int_{\mathbb{D}} |f(z)|^{p} (1-|z|^2)^{\frac{p(2+\beta)}{q}-2}\, dA(z),
\end{align*}
where in the third inequality we have used that
\[
(a+b)^{p} \leq 2^{p} (a^p+b^{p})
\]
for every $a, b \geq 0$.\ Then $f \in A^p_{\frac{p(2+\beta)}{q}-2}$.\ 

We conclude that\[
    A^q_{\beta}=\bigcap_{g\in V_{q,p}}\left[T_{g},A^p_{\alpha}\right]\ (p\le q),\qquad A^p_{\frac{p(2+\beta)}{q}-2}=\bigcap_{g\in V_{q,p}}\left[T_{g},A^p_{\alpha}\right]\ (q<p).
\] If $q<p$, observe that\[
    \frac{2+\beta}{q}= \frac{\frac{p(2+\beta)}{q}-2+2}{p}.
\] Hence, by Lemma~\ref{0}(a), we have $A^q_{\beta} \subseteq A^p_{\frac{p(2+\beta)}{q}-2}$.\ The inclusion is strict, since the reverse inclusion would force $p<q$ by Lemma~\ref{0}\textup{(b)}, a
contradiction.\ Which completes the proof.
\end{proof}

\begin{proof}[\textbf{Proof of Theorem~\ref{theo:3}}]
For $q> p$, the inclusion $A^q_{\beta}\subset [T_g,A^p_{\alpha}]$ is immediate from the boundedness of $T_{g}:A^q_{\beta}\to A^p_{\alpha}$.\ Suppose, for contradiction, that equality holds for some $g$.\ Then the identity map
\[
\operatorname{id}: A^{q}_{\beta}\to [T_g,A^p_{\alpha}]
\]
is a continuous linear bijection between complete quasi-normed spaces.\ By the open mapping theorem (for quasi-Banach spaces if necessary;\ see \cite{19}), it follows that the norms $\|f\|_{\beta,q}$ and $\|T_{g}f\|_{p,\alpha}$ are equivalent for $f\in A^{q}_{\beta}$.\ Consequently $T_{g}:A^q_{\beta}\to A^{p}_{\alpha}$ is bounded below.

However, it is known~\cite[Theorem~1.1]{21} that for this range of indices every bounded operator $T_{g}:A^q_{\beta}\to A^{p}_{\alpha}$ is necessarily compact.\ A compact operator on an infinite‑dimensional quasi-Banach space cannot be bounded below — a contradiction.\ Hence the inclusion is strict. 

For $q<p$, by Theorem~\ref{theo:2}, we obtain $A^q_{\beta} \subsetneq \bigcap_{g\in V_{q,p}}\left[T_{g},A^p_{\alpha}\right]\subseteq [T_{g},A^p_{\alpha}]$.
\end{proof}

To prove Theorem~\ref{theo:5}, we first establish norm inequalities for $T_g:A^p_{\beta}\to A^{p}_{\alpha}$ bounded above and below.
\begin{proposition}
\label{pro}
Suppose $T_g : A_\beta^p\to A_\alpha^p$ is bounded above and below.\ That is, there are constants $0<c_1\le c_2<\infty$ such that
    \begin{equation}
        c_1\int_{\mathbb D}|f(z)|^p\,dA_{\beta}(z)
\le
\int_{\mathbb D}|T_gf(z)|^p\,dA_{\alpha}(z)
\le
c_2\int_{\mathbb D}|f(z)|^p\,dA_{\beta}(z),
\qquad f\in A_\beta^p.\label{eq1}
    \end{equation}
    Then we have 
    \[
        c_1\int_{\mathbb D}|f(z)|^p\,dA_{\beta}(z)
\le
\int_{\mathbb D}|T_gf(z)|^p\,dA_{\alpha}(z)
\le
c_2\int_{\mathbb D}|f(z)|^p\,dA_{\beta}(z),
\qquad f\in \operatorname{Hol}(\mathbb D).
\]
\end{proposition}
The proof of Proposition~\ref{pro} rests on the following result.
\begin{lemma}
\label{lem:1}
Let $0< r< R< 1$, and suppose that $T_g : A_\beta^p\to A_\alpha^p$ satisfies \textup{\eqref{eq1}}.\ Then there exists a constant $C>0$ such that
\begin{equation}
\int_{\Delta(a,r)}|f(z)|^p(1-|z|^2)^{\beta}\,dA(z)
\le
C\int_{\Delta(a,R)}
|f(z)g'(z)|^p(1-|z|^2)^{\alpha+p}\,dA(z),
\label{eq6}
\end{equation}
for every $a\in\mathbb{D}$ and every $f\in\operatorname{Hol}(\mathbb{D})$.
\end{lemma}
\begin{proof}
Fix $a\in\mathbb D$ and $f\in\operatorname{Hol}(\mathbb D)$.\ Applying the change of variables $z=\varphi_a(\zeta)$ and the standard identities\[
1-|\varphi_a(\zeta)|^2
=
\frac{(1-|a|^2)(1-|\zeta|^2)}
     {|1-\overline a\zeta|^2},
\qquad
|\varphi_a'(\zeta)|^2
=
\frac{(1-|a|^2)^2}{|1-\overline a\zeta|^4},
\]
we obtain
\[
\begin{aligned}
\int_{\Delta(a,r)}|f(z)|^p(1-|z|^2)^{\beta}\,dA(z)
&=
(1-|a|^2)^{\beta+2}
\int_{|\zeta|<r}
|f(\varphi_a(\zeta))|^p
\frac{(1-|\zeta|^2)^\beta}
     {|1-\overline a\zeta|^{2\beta+4}}\,dA(\zeta).
\end{aligned}
\]
The last weight factor is uniformly comparable to $1$ on
$|\zeta|<r$.\ Therefore,
\[
\int_{\Delta(a,r)}|f(z)|^p(1-|z|^2)^{\beta}\,dA(z)
\sim
(1-|a|^2)^{\beta+2}
\int_{|\zeta|<r}|f(\varphi_a(\zeta))|^p\,dA(\zeta).
\]
Similarly, arguing as above, we have
\[
\begin{aligned}
\int_{\Delta(a,R)}
|f(z)g'(z)|^p(1-|z|^2)^{\alpha+p}\,dA(z)\sim
(1-|a|^2)^{p+2+\alpha}
\int_{|\zeta|<R}
|f(\varphi_a(\zeta))g'(\varphi_a(\zeta))|^p
\,dA(\zeta).
\end{aligned}
\]The comparison constants are independent of $a$ and $f$.\ Now we set \[
H_a(\zeta)
=(1-|a|^2)^{1+\frac{\alpha-\beta}{p}} g'(\varphi_a(\zeta))\quad \zeta \in \mathbb{D}.
\]Thus the proof of~\eqref{eq6} reduces to showing that there exists a constant $C_1>0$ such that \begin{equation}
\int_{|\zeta|<r}|h(\zeta)|^p\,dA(\zeta)
\le
C_1\int_{|\zeta|<R}
|h(\zeta)H_a(\zeta)|^p\,dA(\zeta)
\label{eq2}
\end{equation}for every $a\in\mathbb{D}$ and every $h\in\operatorname{Hol}(\mathbb{D})$.

We first claim that $H_a$ are uniformly bounded on compact subsets of $\mathbb D$ and that there exist $0<r_0<1$ and $c_0>0$ such that for each $a$ we can choose $\zeta_a$ with $|\zeta_a|<r_0$ and
$|H_a(\zeta_a)|\ge c_0$.\ 
 
Let $u(z)=(1-|z|^2)^{1+\frac{\alpha-\beta}{p}} |g'(z)|$.\ Substituting $z=\varphi_a(\zeta)$ into the definition of $u$ and using the automorphism identity, we derive
\[
|H_a(\zeta)|
=
u(\varphi_a(\zeta))
\left(
\frac{|1-\overline a\zeta|^2}{1-|\zeta|^2}
\right)^{1+\frac{\alpha-\beta}{p}}.
\] Since $T_g$ is bounded, we have $u\in L^\infty$ and $1+\frac{\alpha-\beta}{p}\ge 0$;\ set $
M=\sup_{z\in\mathbb D}u(z)<\infty$.\ For every fixed $0<t<1$, $\left(\frac{|1-\overline a\zeta|^2}{1-|\zeta|^2}
\right)^{1+\frac{\alpha-\beta}{p}}$ is bounded uniformly for $a\in\mathbb D$ and
$|\zeta|\le t$.\ Therefore $\{H_a:a\in\mathbb D\}$ is uniformly bounded on compact subsets of $\mathbb D$.

To produce the point $\zeta_a$, we consider the normalized kernel
\[
k_a(z)=\frac{(1-|a|^2)^{(\beta+2)/p}}{(1-\overline a z)^{2(\beta+2)/p}},
\qquad z\in \mathbb{D},
\]so that $\|k_a\|_{\beta,p}=1$.\ Making the change of variables, we obtain, for every $0<t<1$,
\begin{equation}
\begin{aligned}
\int_{\mathbb D\setminus\Delta(a,t)}
|k_a(z)|^p\,dA_\beta(z)
&=\int_{t\le|\zeta|<1}|k_a(\varphi_a(\zeta))|^p
(\beta+1)(1-|\varphi_a(\zeta)|^2)^\beta
|\varphi_a'(\zeta)|^2\,dA(\zeta)\\
&=\int_{t\le|\zeta|<1}dA_\beta(\zeta)\\
&=(1-t^2)^{\beta+1}.
\end{aligned}\label{eq3}
\end{equation} Thus the integral tends to zero as $t\to1$, uniformly in $a$.

By the derivative characterization and the definition of $u$, there exists $c_3>0$ such that
\[
\int_{\mathbb D}|T_gf|^p\,dA_\alpha
\le c_3\int_{\mathbb D}|uf|^p\,dA_\beta,
\qquad f\in A_\beta^p.
\]
We can choose $r_0\in(0,1)$ sufficiently close to $1$ so that $
M^p(1-r_0^2)^{\beta+1}<\frac{c_1}{2c_3}$.\ Substituting the test function $k_a$ into~\eqref{eq1}, and using~\eqref{eq3}, we have
\[
\begin{aligned}
\int_{\Delta(a,r_0)}
|k_a(z)|^p u(z)^p\,dA_\beta(z)
\ge
\frac{c_1}{c_3}-
\int_{\mathbb D\setminus\Delta(a,r_0)}
|k_a(z)|^p u(z)^p\,dA_\beta(z) >\frac{c_1}{2c_3}.
\end{aligned}
\]Since $\int_{\Delta(a,r_0)}|k_a(z)|^p\,dA_\beta(z)\le1$, arguing by contradiction we see that there exists $z_a\in\Delta(a,r_0)$ such that
\[
u(z_a)> \left(\frac{c_1}{4c_3}\right)^{1/p}.
\]

Furthermore, let $\zeta_a=\varphi_a(z_a)$.\ Then $|\zeta_a|<r_0$, and the preceding comparison gives a constant $c_0>0$ such that
\[
|H_a(\zeta_a)|\ge c_0,
\qquad a\in\mathbb D.
\]

Next, we prove \eqref{eq2}.\ Suppose that~\eqref{eq2} fails.\ After normalization, there are sequences
$\{a_n\}_{n\ge 1}\subset\mathbb D$ and $\{h_n\}_{n\ge 1}\subset\operatorname{Hol}(\mathbb D)$ such that
\begin{equation}
\int_{|\zeta|<r}|h_n(\zeta)|^p\,dA(\zeta)=1\ (n\ge 1),
\qquad
\int_{|\zeta|<R}
|h_n(\zeta)H_{a_n}(\zeta)|^p\,dA(\zeta)\to0\ (n\to \infty).
\label{eq5}
\end{equation}
By Montel's theorem, after passing to a subsequence we may assume that $H_{a_n}\to H$ uniformly on compact subsets of $\mathbb D$.\ For each $a_n$, choose $\zeta_n$ with $|\zeta_n|<r_0$ satisfying $|H_{a_n}(\zeta_n)|\ge c_0$.\ Passing to a further subsequence, we may assume that
\[
\zeta_n\to\zeta_*,
\qquad |\zeta_*|\le r_0.
\]
By the local uniform convergence and the choice of $\zeta_n$, $|H(\zeta_*)|\ge c_0$, hence $H\not\equiv0$.\ Since the zeros of $H$ are isolated, we can choose $s\in(r,R)$ such that $H$ has no zero on the circle $|\zeta|=s$.\ We set $
m_s=\min_{|\zeta|=s}|H(\zeta)|>0$.\ For all sufficiently large $n$ and all $\zeta$ with $|\zeta|=s$, uniform convergence on this circle gives $|H_{a_n}(\zeta)|\ge m_s/2$.

Since $s<R$, subharmonicity of $|h_nH_{a_n}|^p$ implies that
\[
\sup_{|\zeta|= s}|h_n(\zeta)H_{a_n}(\zeta)|^p
\le
C_{s,R}
\int_{|w|<R}|h_n(w)H_{a_n}(w)|^p\,dA(w)
\to0.
\]
Using the lower bound for $|H_{a_n}(\zeta)|$ on $|\zeta|=s$, we obtain
\[
\sup_{|\zeta|=s}|h_n(\zeta)|
\le
\frac{2}{m_s}\sup_{|\zeta|=s}|h_n(\zeta)H_{a_n}(\zeta)|
\to0.
\]
The maximum modulus principle now shows that $h_n\to0$ uniformly on $|\zeta|\le r$, contradicting the normalization in~\eqref{eq5}.\ This proves~\eqref{eq2} and hence~\eqref{eq6}.
\end{proof}
\begin{proof}[\textbf{Proof of Proposition~\ref{pro}}]
Let $f\in\operatorname{Hol}(\mathbb D)$.\ For $0<r<R<1$, we can choose a countable $r$-lattice $\{b_j\}_{j\ge1}$ such that
\[
\mathbb D=\bigcup_j\Delta(b_j,r),
\qquad
\sum_j\mathbf 1_{\Delta(b_j,R)}(z)\le N,
\quad z\in\mathbb D,
\]
where $N<\infty$ is independent of $z$ (see \cite[page~68, Lemma~15]{10}) and \cite[Lemma 4.8]{24}).\ Using the covering property, summing~\eqref{eq6} over $j$, then applying the finite-overlap bound, we obtain
\begin{equation}
\begin{aligned}
\int_{\mathbb D}|f(z)|^p(1-|z|^2)^{\beta}\,dA(z)
&\le
\sum_j
\int_{\Delta(b_j,r)}|f(z)|^p(1-|z|^2)^{\beta}\,dA(z)\\
&\le
C\sum_j
\int_{\Delta(b_j,R)}
|f(z)g'(z)|^p(1-|z|^2)^{\alpha+p}\,dA(z)\\
&= C\int_{\mathbb D}
|f(z)g'(z)|^p(1-|z|^2)^{\alpha+p}
\sum_j\mathbf 1_{\Delta(b_j,R)}(z)\,dA(z)\\
&\le
CN\int_{\mathbb D}
|f(z)g'(z)|^p(1-|z|^2)^{\alpha+p}\,dA(z).
\end{aligned}
\label{eq8}
\end{equation}

If $f\in A_\beta^p$, we already have both inequalities in \eqref{eq1} by assumption.\ Suppose that $f\notin A_\beta^p$.\ If $T_gf\in A_\alpha^p$, then the derivative characterization and \eqref{eq8} give $f\in A_\beta^p$, a contradiction.\ Thus both integrals in \eqref{eq1} are infinite in this case, so the inequalities hold for every
$f\in\operatorname{Hol}(\mathbb D)$.
\end{proof}
Proposition~\ref{pro} enables us to prove that it is possible that $A^p_{\beta}=[T_{g},A^p_{\alpha}]$.
\begin{proof}[\textbf{Proof of Theorem~\ref{theo:5}}]
Since $T_g$ is bounded, it follows that $A^p_{\beta}\subset[T_{g},A^p_{\alpha}]$.\ When  $T_g$ is bounded below, by Proposition~\ref{pro} we have $A^p_{\beta}=[T_{g},A^p_{\alpha}]$.

If $T_g$ is not bounded below, assume that $A^p_{\beta}=[T_{g},A^p_{\alpha}]$.\ Then we obtain\[
    \|f\|_{p,\beta}\sim \|f\|_{[T_{g},A^p_{\alpha}]}
\] by the open mapping theorem \cite{19}.\ Therefore, this contradicts the fact that $T_g$ is not bounded below.
\end{proof}

The space $[T_g, A^p_{\alpha}]$ can be regarded as a weighted Bergman space with a non-radial weight.\ To fix notation, let $w:\mathbb{D}\to[0,\infty)$ be a weight function.\ The associated weighted Bergman space is
\[
A_w^p = \biggl\{ f\in\operatorname{Hol}(\mathbb{D}) : \|f\|_{A_w^p}^p := \int_{\mathbb{D}} |f(z)|^p\,w(z)\,dA(z) < \infty \biggr\}.
\] 
Now, given $g\in V_{q,p}$, take the concrete weight $w(z)=|g'(z)|^p(1-|z|^2)^{\alpha+p}$.\ Then
\[
[T_g, A^p_{\alpha}] = A_w^p.
\]

It is well known \cite{24} that polynomials are dense in the weighted Bergman space $A^p_\alpha$.\ However, things are more complicated for $[T_g, A^p_{\alpha}]$ when considering all $g\in V_{q,p}$.\ In Proposition~\ref{pro:1}, we give an equivalent condition for the polynomials are dense in $[T_g, A^p_{\alpha}]$.\ In Proposition~\ref{pro:2}, we are given an example showing that polynomials are not dense in $ [T_{g},A^p_{\alpha}] $.\ Fortunately, by an argument similar to that in \cite{6} and exploiting the properties of meromorphic optimal domains, we can finally prove that every holomorphic optimal domain $[T_g, A^p_{\alpha}]$ is separable.\ The proof is straightforward and is omitted.

\begin{proposition}
\label{pro:1}
The polynomials are dense in $[T_g, A^p_{\alpha}]$ if and only if $A^q_{\beta}$ is dense in $[T_g, A^p_{\alpha}]$.
\end{proposition}

\begin{proof}[\textbf{Proof}]
First, if the polynomials are dense in $[T_g,A^p_{\alpha}]$, then $A^q_{\beta}$ is certainly dense in $[T_g,A^p_{\alpha}]$.

Conversely, suppose that $A^q_{\beta}$ is dense in $[T_g,A^p_{\alpha}]$.\ For any $f \in [T_g,A^p_{\alpha}]$ and $\varepsilon>0$, choose $h\in A^q_{\beta}$ with $\|f-h\|_{[T_g,A^p_{\alpha}]}<\varepsilon/2$.\ Since the polynomials are dense in $A^q_{\beta}$, we may choose a polynomial $P$ satisfying $\|h-P\|_{q,\beta}<\varepsilon/(2C)$, where $C=\|T_g\|_{A^q_\beta\to A^p_\alpha}$.\ Thus 
\[\|f - P\|_{[T_g,A^p_{\alpha}]}\lesssim \|f - h\|_{[T_g,A^p_{\alpha}]} +C\|h-P\|_{q,\beta}<\varepsilon,\]\ so the polynomials are dense in $[T_g,A^p_{\alpha}]$.
\end{proof}
Non‑radial weighted Bergman spaces remain largely unexplored;\ only a handful of results are known (see, e.g., \cite{10}).\ Even so, the classical theory of weighted Bergman spaces $A_\alpha^p$ allows us to construct an explicit symbol $g\in V_{q,p}$ for which polynomials are not dense in $[T_g,A^p_\alpha]$, as demonstrated in \cite[Proposition 4]{7}.

We recall that a function $f\in A^p_{\alpha+p}$ is called cyclic
if there is a sequence of polynomials $\{p_n\}_{n\ge1}$ such that $fp_n\to1$ in $A^p_{\alpha+p}$.\ By \cite[Theorem~1.1]{16}, there exists a singular inner function
that is not cyclic in $A^p_{\alpha+p}$.\
This leads to the following proposition.
\begin{proposition}
\label{pro:2}
Suppose that $g(z) = \int_0^z S(\zeta)d\zeta$, where $S$ is a non-cyclic singular inner function in $A_{\alpha+p}^p$.\ Then the polynomials are not dense in $[T_g, A^p_{\alpha}]$.
\end{proposition}
\begin{proof}[\textbf{Proof}]
Let $g(z) = \int_0^z S(\zeta)\,d\zeta$.\ It is readily verified that $g \in V_{q,p}$ and $\frac{1}{S}$ belongs to $[T_g,A^p_{\alpha}]$.\ We next prove that $\frac{1}{S}$ cannot be approximated by polynomials in the norm of $[T_g,A^p_{\alpha}]$.\ Suppose, to the contrary, that there exists a sequence of polynomials $\{p_n\}_{n\ge 1}$ such that
\[
\left\|\frac{1}{S} - p_n\right\|_{[T_g,A^p_{\alpha}]} \to 0 \quad (n\to\infty).
\]
Using the same norm equivalence, this implies $\|1 - p_n S\|_{p,\alpha+p} \to 0$.\ Thus $S$ is cyclic in $A_{\alpha+p}^p$, which contradicts the choice of $S$.
\end{proof}
\section{On pointwise multipliers and inclusion relations of \texorpdfstring{$[T_g,A^p_{\alpha}]$}{[Tg,Ap alpha]}}
\quad Having proved in Theorem~\ref{theo:1} that $\big([T_{g},A^p_{\alpha}],\|\cdot\|_{[T_{g},A^p_{\alpha}]}\big)$ is a quasi-Banach space, a natural objective is to characterize its pointwise multipliers.\ In this section, we prove that if $g\in V_{q,p_2}$ and $p_1<p_2$, $[T_g,A^{p_2}_{\alpha}]$ is a proper subspace of $[T_g, A^{p_1}_{\alpha}]$, and we establish the parameter characterization of the embedding $[T_g, A^{s}_{\gamma}] \subset [T_g, A^{p}_{\alpha}]$ for $g\in \operatorname{Hol}(\mathbb D)$ with $g'$ having finitely many zeros.\ The multiplier spaces are of critical importance in the proof of strict containment.
\begin{proof}[\textbf{Proof of Theorem~\ref{theo:6}}]
Let $f\in [T_g,A^p_{\alpha}]$ and $h\in H^\infty$.\ A direct computation yields
\begin{align*}
 \|T_g(hf)\|_{p,\alpha} 
\le C \|T_g(f)\|_{p,\alpha} < \infty,
\end{align*}
with $C \sim \|h\|_\infty$.\ Thus $M_h$ is bounded on $[T_g,A^p_{\alpha}]$ for every $h\in H^\infty$.

Conversely, suppose $h\in \operatorname{Hol}(\mathbb D)$ and $M_h$ is bounded on $[T_g,A^p_{\alpha}]$.\ For a fixed $z\in\mathbb D$, let $\delta_z$ denote the point evaluation functional on $[T_g,A^p_{\alpha}]$.\ By the proof of Theorem~\ref{theo:1}, $\delta_z$ is bounded and, since the space contains the constant function, $\|\delta_z\|\neq 0$.\ For each $f\in [T_g,A^p_{\alpha}]$ with $\|f\|_{[T_g,A^p_{\alpha}]}=1$, we have \[
    |h(z)f(z)| = |\delta_z(M_h f)| \le \|\delta_z\|\, \|M_h\|.
\] Taking the supremum over all such normalized $f$, we get
\[
|h(z)|\, \|\delta_z\|
= \sup_{\|f\|=1} |h(z)f(z)|
\le \|M_h\|\, \|\delta_z\|.
\]
 It follows that $|h(z)| \le \|M_h\|$ for each $z\in\mathbb D$.\ Hence $h\in H^\infty$.
\end{proof}

The following lemma is an extension of \cite[Lemma 2]{8}, where only the case
$0<\kappa<1$ with $0<\varepsilon<\min\{\kappa,1-\kappa\}$ was treated.
\begin{lemma}\label{lem:2}
Let $\kappa\ge0$, $\varepsilon> 0$, and let $g\in \mathcal{A}^{-\kappa}$ with
$g\not\equiv 0$.\ Then there exists $\varphi\in\operatorname{Hol}(\mathbb D)$
such that
\[
\varphi g\in \mathcal{A}^{-(\kappa+\varepsilon)}\setminus \mathcal{A}^{-\kappa}.
\]
\end{lemma}

\begin{proof}
By \cite[Theorem~1.2]{1}, we may choose
$f_1,f_2\in \mathcal{A}^{-\varepsilon}$ such that
\[
|f_1(z)|+|f_2(z)|\sim(1-|z|^2)^{-\varepsilon}.
\]
Fix an integer $N\ge1$ with $N\varepsilon>\kappa$.\ We now prove that $gf_i^N\notin \mathcal{A}^{-\kappa}$ for $i=1$ or $2$.\ If both $gf_1^N$ and $gf_2^N$ belong to $\mathcal{A}^{-\kappa}$, then
\[
|g(z)|(1-|z|^2)^{-N\varepsilon}
\lesssim
|g(z)|\bigl(|f_1(z)|^N+|f_2(z)|^N\bigr)
\lesssim(1-|z|^2)^{-\kappa},
\]
where we used
$(a+b)^{N} \leq 2^{N} (a^N+b^{N})$ for $a,b\ge0$.\ Thus
\[
|g(z)|\lesssim(1-|z|^2)^{N\varepsilon-\kappa},
\]
which forces $g\equiv0$ by the maximum modulus principle, a contradiction.

Hence we may assume that $gf_1^N\notin \mathcal{A}^{-\kappa}$, and set $
m=\min\{1\le k\le N:gf_1^k\notin \mathcal{A}^{-\kappa}\}$.\ Since $g\in \mathcal{A}^{-\kappa}$, we have $gf_1^{m-1}\in \mathcal{A}^{-\kappa}$.\ Consequently,
\[
gf_1^m=(gf_1^{m-1})f_1
\in \mathcal{A}^{-(\kappa+\varepsilon)}\setminus \mathcal{A}^{-\kappa}.
\]
Taking $\varphi=f_1^m$ proves the assertion for all
$\kappa\ge0$ and $\varepsilon>0$.
\end{proof}
\begin{proof}[\textbf{Proof of Theorem~\ref{theo:7}}]
Since $p_2>p_1>0$, the classical inclusion $A^{p_2}_{\alpha}\subset A^{p_1}_{\alpha}$ holds.\ The definition of the optimal domain yields the containment
\[
[T_g,A^{p_2}_{\alpha}]\subset [T_g,A^{p_1}_{\alpha}].
\]

Suppose to the contrary that equality holds, i.e. $[T_g,A^{p_2}_{\alpha}]=[T_g,A^{p_1}_{\alpha}]$, and set $\kappa=1 + \frac{2+\alpha}{p_2} - \frac{2+\beta}{q}$.\ Then by the assumption of the theorem, $\kappa\ge 0$.

\textbf{Case 1:} $q \le p_1 < p_2$.\ We consider the multiplication operator $M_h$ acting from $A_\beta^q$ to $[T_g,A_\alpha^{p_2}]$.\ Using the definition of $T_g$ and the derivative property of Bergman spaces, it follows that $g'h$ is a multiplier from $A^q_\beta$ to $A^{p_2}_{\alpha+p_2}$.\ By \cite[Theorem~1]{23}, the multiplier space $X_1$ from $A^q_\beta$ to $[T_g, A_\alpha^{p_2}]$ is given by
\[
\left\{ h \in \operatorname{Hol}(\mathbb{D}) : \sup_{z \in \mathbb{D}} (1-|z|^2)^{\kappa} |g'(z)h(z)| < \infty \right\},
\]
and the multiplier space $X_2$ from $A_\beta^q$ to $[T_g, A_\alpha^{p_1}]$ is given by
\[
 \left\{ h \in \operatorname{Hol}(\mathbb{D}) : \sup_{z \in \mathbb{D}} (1-|z|^2)^{1 + \frac{2+\alpha}{p_1} - \frac{2+\beta}{q}} |g'(z)h(z)| < \infty \right\}.
\]

Assume $[T_g,A_\alpha^{p_2}] = [T_g,A_\alpha^{p_1}]$.\ Then the corresponding multiplier spaces must be equal, i.e., $X_1 = X_2$, by the closed graph theorem.\ Since $g \in V_{q,p_2}$, we have 
\[
\sup_{z \in \mathbb{D}} (1-|z|^2)^{\kappa} |g'(z)| < \infty.\] 
Hence the space $X_1$ contains nonzero elements.\ Choose $
\varepsilon = \frac{2+\alpha}{p_1} - \frac{2+\alpha}{p_2}>0$.\ For each fixed $0\not\equiv h \in X_1$, by Lemma~\ref{lem:2} there exists $\varphi \in \operatorname{Hol}(\mathbb{D})$ such that
\[
g' h \varphi \in \mathcal{A}^{-\left(\kappa+\varepsilon\right)} \setminus \mathcal{A}^{-\kappa}.
\]
It then follows that $\varphi h \in X_2 \setminus X_1$, which contradicts that $X_1 = X_2$.

\textbf{Case 2:} $p_1 < q \le p_2$.\ $g$ is the same as in Case~1.\ Similarly, by \cite[Theorem~1]{23}, we obtain the multiplier space $X_1$ from $A_\beta^q$ to $[T_g, A_\alpha^{p_2}]$ is given by
\[
 \left\{ h \in \operatorname{Hol}(\mathbb{D}) : \sup_{z \in \mathbb{D}} (1-|z|^2)^{\kappa} |g'(z)h(z)| < \infty \right\},
\]
and for the multiplier space $X_2$ from $A_\beta^q$ to $[T_g,A_\alpha^{p_1}]$, we have  
\[
X_2 = \left\{ h \in \operatorname{Hol}(\mathbb{D}) : g'h\in A^r_\delta, r=\frac{qp_1}{q-p_1}, \delta=\frac{q(\alpha+p_1)-p_1\beta}{q-p_1}\right\}.
\]

Suppose that $[T_g,A^{p_2}_{\alpha}]=[T_g,A^{p_1}_{\alpha}]$.\ By the argument
of Case~1 we then have $X_1=X_2$.\ Since $g \in V_{q,p_2}$, the space $X_1$ contains nonzero elements.\ From $\frac{2+\alpha}{p_2} - \frac{1}{q}<\frac{1+\alpha}{p_1}$, it follows that we can choose $\varepsilon>0$ with
\[\frac{2+\alpha}{p_2} - \frac{1}{q}+\varepsilon<\frac{1+\alpha}{p_1}
 .
\] For each fixed $0\not\equiv h \in X_1$, by Lemma~\ref{lem:2} there exists $\varphi \in \operatorname{Hol}(\mathbb{D})$ such that
\[
g' h \varphi \in \mathcal{A}^{-\left(\kappa +\varepsilon\right)} \setminus \mathcal{A}^{-\kappa}.
\]
Since $\kappa + \varepsilon<\frac{1+\delta}{r}=1 + \frac{1+\alpha}{p_1} - \frac{1+\beta}{q}$, we obtain $\mathcal{A}^{-\left(\kappa +\varepsilon\right)}\subset A^r_\delta$.\ It then follows that $\varphi h \in X_2 \setminus X_1$, which contradicts that $X_1 = X_2$.

\textbf{Case 3:} $p_1 < p_2 < q$.\ Analogously, by \cite[Theorem~1]{23}, the multiplier space $X_1$ from $A_\beta^q$ to $[T_g, A_\alpha^{p_2}]$ is given by
\[
\left\{ h \in \operatorname{Hol}(\mathbb{D}) : g'h\in A^{r_1}_{\delta_1}, r_1=\frac{qp_2}{q-p_2}, \delta_1=\frac{q(\alpha+p_2)-p_2\beta}{q-p_2}\right\},
\]
and the multiplier space $X_2$ from $A_\beta^q$ to $[T_g, A_\alpha^{p_1}]$ is given by
\[
\left\{ h \in \operatorname{Hol}(\mathbb{D}) : g'h\in A^{r_2}_{\delta_2}, r_2=\frac{qp_1}{q-p_1}, \delta_2=\frac{q(\alpha+p_1)-p_1\beta}{q-p_1}\right\}.
\]

Assume $[T_g, A_\alpha^{p_2}] = [T_g, A_\alpha^{p_1}]$.\ By the argument
of Case~1 we have $X_1=X_2$.\ Since $g \in V_{q,p_2}$, we get $g'\in A^{r_1}_{\delta_1}$.\ Hence $X_1$ contains nonzero elements.\ From $\frac{2+\alpha}{p_2} - \frac{1}{q}<\frac{1+\alpha}{p_1}$, we can choose $\varepsilon>0$ satisfying
\[
 \frac{2+\alpha}{p_2} - \frac{1}{q}+\varepsilon<\frac{1+\alpha}{p_1}.
\]
For each fixed $0\not\equiv h \in X_1$, by Lemma~\ref{lem:2} again, there exists $\varphi \in \operatorname{Hol}(\mathbb{D})$ such that
\[
g' h \varphi \in \mathcal{A}^{-\left(\kappa+\varepsilon\right)} \setminus \mathcal{A}^{-\kappa}.
\]
Since \[\kappa=\frac{2+\delta_1}{r_1},\qquad
\kappa+ \varepsilon<\frac{1+\delta_2}{r_2}=1 + \frac{1+\alpha}{p_1} - \frac{1+\beta}{q},\]
we obtain $A^{r_1}_{\delta_1}\subset \mathcal{A}^{-\kappa}$ and $\mathcal{A}^{-\left(\kappa +\varepsilon\right)}\subset A^{r_2}_{\delta_2}$.\ It then follows that $\varphi h \in X_2 \setminus X_1$, which contradicts that $X_1 = X_2$.\ This completes the proof of conclusions~\textup{(a)} and~\textup{(b)}.

Regarding compactness, we consider a sequence $\{f_n\}_{n\ge 1} \subset [T_g, A_\alpha^{p_2}]$ such that $\|f_n\|_{[T_g, A_\alpha^{p_2}]} \le 1$, consequently, $\|T_g(f_n)\|_{p_2,\alpha} \le 1$.\ Since the embedding $A_{\alpha}^{p_2} \subset A_{\alpha}^{p_1}$ is compact by Lemma~\ref{0}, $T_g(f_n)$ admits a convergent subsequence in $A_{\alpha}^{p_1}$, this is equivalent to $f_n$ possessing a convergent subsequence in $[T_g, A_{\alpha}^{p_1}]$.\ Thus, the embedding is compact, which proves \textup{(c)}. \end{proof}
To prove Theorem~\ref{theo:8}, we first give a proposition.\ For non-constant $g\in\operatorname{Hol}(\mathbb D)$, arguing as in the proof of Theorem~\ref{theo:1} shows that $[T_g,A^p_\alpha]$ is a quasi-Banach space.
\begin{proposition}
\label{pro:3}
Let $0< s<\infty$, $-1<\gamma<\infty$, and let $g \in\operatorname{Hol}(\mathbb{D})$ satisfy $g'\neq0$ on $\mathbb D$.\ If the inclusion $[T_g,A^s_\gamma] \subset[T_g,A^p_\alpha]$ is compact, then the embedding $ A^s_{s+\gamma}\subset A^p_{p+\alpha}$ is compact.
\end{proposition}

\begin{proof}[\textbf{Proof}]
Let $\{h_n\}_{n\ge 1}\subset A^s_{s+\gamma}$ satisfy $\sup_{n}\|h_n\|_{s,s+\gamma}<\infty$.\ Since $g'$ is zero-free on $\mathbb D$, the
functions $f_n=\frac{h_n}{g'}$ are holomorphic in $\mathbb D$.\ We have $f_n\in [T_g,A^s_\gamma]$ and
\[
\|f_n\|_{[T_g,A^s_\gamma]}
=\|T_gf_n\|_{s,\gamma}
\sim \|h_n\|_{s,s+\gamma},
\]
so $\{f_n\}$ is bounded in $[T_g,A^s_\gamma]$.\ By the compactness of the inclusion $[T_g,A^s_\gamma] \subset[T_g,A^p_\alpha]$, there exists a
subsequence $\{f_{n_k}\}_{k\ge 1}$ and a function $f\in[T_g,A^p_\alpha]$ such that \[
\|f_{n_k}-f\|_{[T_g,A^p_\alpha]}
=\|T_g(f_{n_k}-f)\|_{p,\alpha}
\to 0\quad (k\to \infty).
\]
Then
\[
\|h_{n_k}-fg'\|_{p,\alpha+p}
=\bigl\|(f_{n_k}-f)\,g'\bigr\|_{p,\alpha+p}
\sim \bigl\|T_g(f_{n_k}-f)\bigr\|_{p,\alpha}
\to 0 .
\]
Thus $h_{n_k}$ converges to $fg'$ in $A^p_{p+\alpha}$.\ Therefore the
embedding $A^s_{s+\gamma}\subset A^p_{p+\alpha}$ is compact.
\end{proof}
\begin{proof}[\textbf{Proof of Theorem~\ref{theo:8}}]
First, we prove that $[T_g,A^s_\gamma] \subset [T_g,A^p_\alpha]$ if and only if $A^s_\gamma \subset A^p_\alpha$.\ The sufficiency follows immediately from the definition of the optimal domain.\ To prove the necessity, we define the range of $T_g$ as 
\[
\mathcal{R}_g = \bigl\{ F\in \operatorname{Hol}(\mathbb{D}) : F(0)=0,\; F'/g'\in \operatorname{Hol}(\mathbb{D}) \bigr\}.\] 
Because $T_g$ is injective, the map $f\mapsto F=T_g f$ is a bijection from $\operatorname{Hol}(\mathbb{D})$ onto $\mathcal{R}_g$.\ For $f\in\operatorname{Hol}(\mathbb D)$, we then have\ $f\in [T_g,A^s_\gamma]$ if and only if $T_g f \in \mathcal{R}_g\cap A^s_\gamma$, and analogously for $[T_g,A^p_\alpha]$.\ Hence the inclusion $[T_g,A^s_\gamma] \subset [T_g,A^p_\alpha]$ is equivalent to
\begin{equation}\label{eq9}
\mathcal{R}_g\cap A^s_\gamma \subset \mathcal{R}_g\cap A^p_\alpha.
\end{equation}

Suppose, for a contradiction, that \eqref{eq9} holds but $A^s_\gamma\not\subset A^p_\alpha$.\ Pick $F_0\in A^s_\gamma\setminus A^p_\alpha$.\ Let $a_1,\dots,a_m$ be the distinct zeros of $g'$ with multiplicities $n_1,\dots,n_m$.\ Then
\[
\mathcal{R}_g\cap A^s_\gamma = \{ F\in A^s_\gamma : F(0)=0,\; F^{(k)}(a_j)=0 \text{ for } 1\le k\le n_j, 1\le j\le m \},
\]
which is a closed subspace of $A^s_\gamma$.\ By Hermite interpolation, we can choose a polynomial $P$ such that
\[
P(0)=F_0(0),\qquad P^{(k)}(a_j)=F_0^{(k)}(a_j)\quad(1\le k\le n_j, 1\le j\le m).
\]
Set $F=F_0-P$.\ Then $F(0)=0$ and $F^{(k)}(a_j)=0$ for all required indices, so $F\in\mathcal{R}_g$.\ Moreover, $F\in A^s_\gamma$ because $F_0\in A^s_\gamma$ and $P$ is a polynomial.\ On the other hand, since $P\in A^p_\alpha$ whereas $F_0\notin A^p_\alpha$, we conclude that $F\notin A^p_\alpha$.\ Therefore, \[ F\in(\mathcal{R}_g\cap A^s_\gamma)\setminus (\mathcal{R}_g\cap A^p_\alpha),\]
\ contradicting~\eqref{eq9}.\ Hence $A^s_\gamma\subset A^p_\alpha$ is necessary.\ Conclusions~\textup{(a)} and~\textup{(b)} follow from Lemma~\ref{0}.

 Now we prove conclusion~\textup{(c)}, if the strict inequality holds, then the embedding $A^s_\gamma \subset A^p_\alpha$ is compact and $[T_g,A^s_\gamma] \subset [T_g,A^p_\alpha]$ is also compact by the same argument as in Theorem~\ref{theo:7}\textup{(c)}.\

 On the other hand, suppose the embedding $[T_g,A^s_\gamma] \subset [T_g,A^p_\alpha]$ is compact.\ Since $g'$ admits only finitely many zeros in $\mathbb{D}$, there exists a finite Blaschke product $B$ such that $g' = B g_1'$, where $g'_1\in \operatorname{Hol}(\mathbb{D})$ is zero-free.\ By Lemma~\ref{lem:5} in the next section, $[T_{g}, A_{\alpha}^p] = [T_{g_1}, A_{\alpha}^p]$.\ Hence, the embedding $[T_{g_1}, A_{\gamma}^s]\subset [T_{g_1}, A_{\alpha}^p]$ is also compact.\ By Proposition~\ref{pro:3}, this yields the compactness of the embedding $A_{\gamma+s}^s \subset A_{\alpha+p}^p$, and the strict inequality follows from~ Lemma~\ref{0}. \end{proof}
\section{The space \texorpdfstring{$V_g^p$}{Vg p} and equivalence of holomorphic optimal domains}
\quad In this section, we examine the space $V_g^p$ and determine the conditions ensuring the coincidence of the holomorphic optimal domains for two $V_{q,p}$ functions $f$ and $g$.

The complete characterization of $V_g^p$ for general $g\in$ $V_{q,p}$ is currently unknown.\ Nevertheless, we observe that if $g$ satisfies certain conditions, the space $V_g^p$ is independent of the choice of $p$.

Firstly, we shall need the following fact, see \cite{14}.
\begin{lemma}\label{lem:4}
Suppose that $f\in A_\alpha^p$.\ Let $B$ be a Blaschke product whose zero sequence, counted with multiplicities, is a finite union of uniformly separated sequences.\ If $f/B\in\operatorname{Hol}(\mathbb D)$,
then $f/B\in A_\alpha^p$.
\end{lemma}

\begin{proof}[\textbf{Proof of Theorem~\ref{theo:9}}]
By the definition of $V_g^p$, the constant functions belong to $V_g^p$.\ Let $f-f(0)=T_{g_1}(h)$ with $h\in H^\infty$.\ Then, for every $z\in\mathbb{D}$, 
$f^{\prime}(z)=h(z)g_1^{\prime}(z)$.\ To establish $f\in V_g^p$, we show $[T_f,A^p_{\alpha}]\supseteq[T_g,A^p_{\alpha}]$.\ Let $u\in [T_g,A^p_{\alpha}]$.\ Then, due to Theorem~\ref{theo:6}, $hu\in[T_g,A^p_{\alpha}]$ and
\begin{align*}
    \|u\|_{[T_f, A^p_{\alpha}]} = \|T_f(u)\|_{p,\alpha} 
    \sim \|f'u\|_{p,\alpha+p}= \|\frac{hg'u}{B}\|_{p,\alpha+p}< \infty,
\end{align*}
 where in the last inequality, we use the fact that $hg'u\in A^p_{\alpha+p}$ and Lemma~\ref{lem:4} , thus $[T_f,A^p_{\alpha}]\supseteq[T_g,A^p_{\alpha}]$, proving the first inclusion.

On the other hand, let $f\in V_g^p$.\ As $g_1^{\prime}$ has no zeros, $
\frac{f^{\prime}}{g_1^{\prime}}\in\operatorname{Hol}(\mathbb{D})$.\ For any $u\in[T_g,A^p_{\alpha}]\subset [T_f,A^p_{\alpha}]$, we know $T_f(u)\in A^p_{\alpha}$ and
\[
\|T_g\left(\frac{f^{\prime}}{g_1^{\prime}}u\right)\|_{p,\alpha}\sim \|Bf'u\|_{p,\alpha+p}\lesssim \|T_f(u)\|_{p,\alpha}<\infty.
\]
Consequently, $f^{\prime}/g_1^{\prime}$ is a multiplier of $[T_g,A^p_{\alpha}]$ and, according to Theorem~\ref{theo:6}, we have $f^{\prime}/g_1^{\prime}\in H^\infty$.\ This proves that $f-f(0)\in T_{g_1}(H^\infty)$, which concludes the theorem.
\end{proof}

Theorem~\ref{theo:9} implies that $V_g^p=T_{g_1}(H^\infty)+\mathbb{C}$, showing that $V_g^p$ is independent of $p$.\ However, we do not have this characterization for all $g\in V_{q,p}$.\ Our next goal is to determine when two holomorphic optimal domains are equal.
\begin{corollary}\label{cor:1}
Suppose that $f,g\in$ $V_{q,p}$ are non-constant and that both satisfy the hypotheses of Theorem~\ref{theo:9}, write $f'=B_1f'_1$ and  $g'=B_2g'_1$.\ Then the following are equivalent:
\begin{enumerate}
\item[\textup{(a)}] $[T_{g},A^p_{\alpha}]=[T_{f},A^p_{\alpha}]$;
\item[\textup{(b)}] There exist two functions $h_1\in ,h_2\in H^\infty$ such that
\[
g=T_{f_1}(h_1)+g(0)\quad\text{and}\quad f=T_{g_1}(h_2)+f(0).
\]Moreover we have that $h_1h_2=B_1B_2$.
\end{enumerate}
\end{corollary}
\begin{proof}[\textbf{Proof}]
If \textup{(a)} holds, then we have $g\in V_f^p$ and $f\in V_g^p$.\ Theorem~\ref{theo:9} therefore yields \textup{(b)}.\ Conversely, \textup{(b)} and Theorem~\ref{theo:9} imply
$g\in V_f^p$ and $f\in V_g^p$, which give the conclusion \textup{(a)}.

Differentiating the identities in \textup{(b)}, we obtain
$g'=h_1f_1'$ and $f'=h_2g_1'$.\ Since $f'=B_1f_1'$ and $g'=B_2g_1'$, we conclude that $h_1h_2=B_1B_2$.
\end{proof}

The following corollary is an immediate consequence of Corollary~\ref{cor:1}.
\begin{corollary}\label{cor:2}
Let $g'=Bh'$ with $h\in V_{q,p}$, $h'\neq0$ on $\mathbb D$, and $B$ a Blaschke product whose zeros, counted with multiplicity, form a finite union of uniformly separated sequences.\ Then we have $[T_g,A^p_{\alpha}]=[T_u,A^p_{\alpha}]$.
\end{corollary}
\begin{proof}[\textbf{Proof}]
By the definition of $g$, we have that $B_1=1$ and $B_2=B$.\ We can choose $k_1=B$ and $k_2=1$, so
    \[
\begin{cases}
g^{\prime}=k_1h^{\prime}\\
h^{\prime}=k_2h^{\prime}
\end{cases}.
\]
By Corollary~\ref{cor:1}, we have the desired conclusion.
\end{proof}

In the remainder of this section we characterize, under auxiliary hypotheses, those symbols $g_1,g_2\in V_{q,p}$ for which the holomorphic optimal domains coincide.
\begin{proposition}
\label{pro:4}
Let $g_1,g_2\in V_{q,p}$ be non-constant such that $g_2^{\prime}-g_1^{\prime}=h^{\prime}\in\operatorname{Hol}(\mathbb{D})$.\ Then $[T_{g_1},A^p_{\alpha}]=[T_{g_2},A^p_{\alpha}]$ if and only if $h\in V^p_{g_1}\cap V^p_{g_2}$.
\end{proposition}
\begin{proof}[\textbf{Proof}]
Suppose $[T_{g_1},A^p_{\alpha}]=[T_{g_2},A^p_{\alpha}]$ and $f\in [T_{g_1},A^p_{\alpha}]$.\ Then
    \[
        \|f\|_{[T_{h},A^p_{\alpha}]}=\left\|T_{h}(f)\right\|_{p,\alpha}\lesssim\|h'f\|_{p,\alpha+p}\lesssim \|g'_2f\|_{p,\alpha+p}+\|g'_1f\|_{p,\alpha+p}<\infty.
    \]
We have that $h\in V^p_{g_1}$, and similarly $h\in V^p_{g_2}$.

Conversely, suppose $h\in V^p_{g_1}\cap V^p_{g_2}$ and $f\in [T_{g_1},A^p_{\alpha}]$, then
\[
    \|f\|_{[T_{g_2},A^p_{\alpha}]}=\left\|T_{g_2}(f)\right\|_{p,\alpha}\lesssim\|g'_2f\|_{p,\alpha+p}\lesssim \|g'_1f\|_{p,\alpha+p}+\|h'f\|_{p,\alpha+p}<\infty.
\]
This implies $[T_{g_1},A^p_{\alpha}]\subset [T_{g_2},A^p_{\alpha}]$.\ By symmetry, $[T_{g_2},A^p_{\alpha}]\subset [T_{g_1},A^p_{\alpha}]$.
\end{proof}

We next consider symbols $g$ whose derivatives are holomorphic on the closed unit disk.\ We need the following auxiliary lemma.

\begin{lemma}
\label{lem:5}
Let $g\in \operatorname{Hol}(\mathbb{D})$.\ If $g^{\prime}=Q u^{\prime}$ with $u\in \operatorname{Hol}(\mathbb{D})$ and $Q$ a holomorphic function on $\bar{\mathbb{D}}$ with no zeros on $\mathbb{T}$, then $[T_g,A^p_{\alpha}]=[T_u,A^p_{\alpha}]$.
\end{lemma}

\begin{proof}
Let $f \in [T_u, A^p_{\alpha}]$.\ Since $Q\in H^\infty$, we have $Qu'f\in A_{\alpha+p}^p$ and hence
$Qf\in[T_u,A_\alpha^p]$.\ Therefore $T_g(f) = T_u(Qf) \in A^p_{\alpha}$, which implies
$[T_u, A^p_{\alpha}] \subset [T_g, A^p_{\alpha}]$.

Now we take $f \in [T_g, A^p_{\alpha}]$, that is, $g'f \in A^p_{\alpha+p}$.\ Since $Q$ has no zeros on $\mathbb{T}$, there exist $r\in(0,1)$ and $\delta>0$ such that $|Q(z)| \ge \delta$ for all $z$ with $r \le |z| \le 1$.\ On the other hand, since $u'f$ is holomorphic on $\mathbb{D}$, there is a constant $C>0$ such that $|u'(z)f(z)| < C$ whenever $|z| \le r$.\ We then have
\begin{align*}
    \|f\|^p_{[T_u, A^p_{\alpha}]} 
    &\lesssim\int_{|z| \le r} |u'(z)f(z)|^p dA_{\alpha+p}(z) +  \frac{1}{\delta^p} \int_{r \le |z| \le 1} |Q(z)u'(z)f(z)|^p dA_{\alpha+p}(z) \\
 &< C+\frac{\|g'f\|^p_{p,\alpha+p}}{\delta^p}.
\end{align*}
Therefore, we obtain the desired conclusion.\ Finally, if $g$ is constant, then $u$ is also constant.\ Clearly, $[T_g,A^p_{\alpha}]=[T_u,A^p_{\alpha}]$
\end{proof}

\begin{proof}[\textbf{Proof of Theorem~\ref{theo:11}}]
First, we note that $[T_z,A^p_{\alpha}]=A^p_{\alpha+p}$ by the derivative characterization of weighted Bergman spaces.\ The remaining part of \textup{(a)} follows from Lemma~\ref{lem:5}.

For part \textup{(b)}, without loss of generality, we consider the case $\lambda_1=1$.\ The general case follows similarly.\ Let $g'(z)=(1-z)h(z)$, where $h\in \operatorname{Hol}(\mathbb{\bar{D}})$ is bounded.\ For each $f\in A^p_{\alpha+p}$, it is clear that $f\in [T_g,A^p_{\alpha}]$.\ Now, to prove the strict inclusion $A^p_{\alpha+p}\subsetneq [T_g,A^p_{\alpha}]$, define a function $f_0$ on $\mathbb D$ by
$$
f_0(z)=\frac{1}{(1-z)^{1+\frac{\alpha+2}{p}}},\qquad z\in\mathbb D.
$$ Then $f_0\notin A^p_{\alpha+p}$.\ We shall show that $f\in [T_g,A^p_\alpha]$.\ Since $
\frac{1}{(1-z)^{\frac{\alpha+2}{p}}}\in A^p_{\alpha+p}$ and $h$ is bounded, we obtain  
\[
    \|f_0\|_{[T_g,A^p_{\alpha}]}\sim \|\frac{g'(z)}{(1-z)^{1+\frac{\alpha+2}{p}}}\|_{p,\alpha+p}=\|\frac{h(z)}{(1-z)^\frac{\alpha+2}{p}}\|_{p,\alpha+p}<\infty.
\]
Thus the proof is complete.
\end{proof}

\section{Examples of \texorpdfstring{$A^p_\alpha\subsetneq[T_g,A^p_{\alpha}]$}{Ap alpha subset [Tg, Ap alpha]}}
\quad In this section, we give some concrete examples in which $A^p_\alpha\subsetneq[T_g,A^p_{\alpha}]$ for $g\in \mathcal{B}$.\ To this end, we need the following two lemmas.
\begin{lemma}
\label{lem:6}
Let $g \in \mathcal{B}_0$ be non-constant.\ Then the operator $T_g : A^p_{\alpha}\to A^p_{\alpha}$ is not bounded below.
\end{lemma}
\begin{proof} 
Since $T_g$ depends only on $g'$, we may assume without loss of generality that $g(0)=0$.\ Assume, for the sake of contradiction, that $T_g$ is bounded below.\ This means there exists a constant $C > 0$ such that
\[
\|T_g(f)\|_{p,\alpha} \geq C \|f\|_{p,\alpha} \quad \text{for all } f \in A^p_{\alpha}.
\]

It is known from \cite[Corollary 4]{3} that  $e^g\in A^p_{\alpha}$ for every $g\in \mathcal{B}_0$.\ In particular, since $ng\in\mathcal{B}_0$ for each positive integer $n$, the function $f_n=e^{ng}$ belongs to $A^{p}_{\alpha}$.\ Moreover, the following norm estimates hold (with implied constants independent of $n$):
\[
\|f_n\|_{p,\alpha} \sim 1+\|ng'e^{n g}\|_{p,\alpha+p} \quad \text{and} \quad \|T_g(f_n)\|_{p,\alpha} \sim \|g'e^{n g}\|_{p,\alpha+p}.
\]
Therefore, we have
\[
\lim_{n \to \infty} \frac{\|T_g(f_n)\|_{p,\alpha}}{\|f_n\|_{p,\alpha}} = 0.
\]
However, if $T_g$ is bounded below, then \[
\lim_{n \to \infty} \frac{\|T_g(f_n)\|_{p,\alpha}}{\|f_n\|_{p,\alpha}} \geq C > 0.
\]
Thus, we arrive at a contradiction, and conclude that $T_g$ cannot be bounded below. 
\end{proof}
\begin{lemma}\label{lem:7}
Let $g\in\mathcal B$ be non-constant and $h\in\operatorname{Hol}(\mathbb D)$.\ When $-1<\lambda\le\alpha+p$, the multiplication operator
$M_h: A^p_\lambda\to [T_g,A^p_\alpha]$ is bounded if and only if
\[
\sup_{z\in\mathbb D}
(1-|z|^2)^{1+\frac{\alpha-\lambda}{p}}|g'(z)h(z)|<\infty.
\]
\end{lemma}
\begin{proof}
It is easy to check that, $M_h:A^p_\lambda\to[T_g,A^p_\alpha]$ is bounded if and only if $M_{g'h}:A^p_\lambda\to A^p_{\alpha+p}$ is bounded.\ The conclusion follows from \cite[Theorem~1]{23}. 
\end{proof}
\begin{proposition}\label{pro:5}
Let $g\in\mathcal B$ be non-constant.\ If either $g\in\mathcal B_0$, or the zero set of $g'$,\
counted with multiplicity, satisfies the Blaschke condition, then $A^p_\alpha$ is strictly contained in
$[T_g,A^p_\alpha]$.
\end{proposition}
\begin{proof}[\textbf{Proof}]
Since $T_g$ is bounded on $A^p_{\alpha}$ for every $g\in\mathcal{B}$, we have $A^p_\alpha\subset[T_g,A^p_{\alpha}]$, and for all~$f\in A^p_{\alpha}$,
\[
\|f\|_{[T_g,A^p_{\alpha}]}=\|T_g(f)\|_{p,\alpha} \le \|T_g\|_{A^p_{\alpha}\to A^p_{\alpha}}\,\|f\|_{p,\alpha}.
\]
Assume, for the sake of contradiction, that $[T_g,A^p_{\alpha}]=A^p_{\alpha}$.\ By the open mapping theorem \cite{19}, the two norms are equivalent, i.e.
\begin{equation}
\|f\|_{p,\alpha} \sim \|f\|_{[T_g,A^p_{\alpha}]} \quad \text{for all } f \in A^p_{\alpha}.\label{eq10}
\end{equation}

First, suppose that $g\in\mathcal{B}_0$.\ By Lemma~\ref{lem:6}, we have shown that $T_g$ fails to be bounded below, which contradicts~\eqref{eq10}.\ Hence $A^p_{\alpha}\subsetneq [T_g,A^p_{\alpha}]$.

On the other hand, it is known from \cite[Theorem~1]{23} that the multiplier space from $A^p_{\alpha+p}$ to $A^p_{\alpha}$ is $\{0\}$.\ By Lemma~\ref{lem:7}, the multiplier space from $A^p_{\alpha+p}$ to $[T_g, A^p_{\alpha}]$ is \[ X=\lbrace h\in \operatorname{Hol}(\mathbb D): g'h\in H^\infty \rbrace.\] 
Since we assume $A^p_{\alpha}= [T_g, A^p_{\alpha}]$, it follows from~\eqref{eq10} that every $h\in X$ is also a multiplier from $A^p_{\alpha+p}$ to $A^p_{\alpha}$.
    
However, when the zero sequence of $g'$ satisfies the Blaschke condition, we can write $g'=Bg_1$ by \cite[Theorem~2.2]{12}, where $B$ is the associated Blaschke product and
$g_1\in\operatorname{Hol}(\mathbb D)$ is zero-free.\ Setting $h = 1/g_1$, we obtain $h \in X$ but $h \not\equiv 0$, which leads to a contradiction. 
\end{proof}
This proposition, together with Theorem~\ref{theo:5}, immediately yields the following corollary.
\begin{corollary}\label{cor:3}
Let $g\in\mathcal B$ be non-constant.\ If the operator $T_g$ is bounded below on $A^p_\alpha$, then $g'$ has infinitely many zeros, and its zero sequence, counted with multiplicity, fails to satisfy the Blaschke condition.
\end{corollary}

\textbf{Acknowledgments.} 
This work was supported by the National Natural Science Foundation of China (No.~12471093), the Guangdong Basic and Applied Basic Research Foundation
(No.~2024A1515010468), and the Li Ka Shing Foundation STU-GTIIT Joint Research Grant (No.~2024LKSFG06).

\def\bysame{J. Xiao}
\bibliographystyle{amsplain}
\def\MR#1{}
\bibliography{refs}
\par\medskip
\noindent
Zengjian Lou: Department of Mathematics, Shantou University,
515063, Shantou, Guangdong, P.R.~China.\\
Email address: \texttt{zjlou@stu.edu.cn}

\par\medskip
\noindent
Jinbang Wen: Department of Mathematics, Shantou University,
515063, Shantou, Guangdong, P.R.~China.\\
Email address: \texttt{23jbwen@stu.edu.cn}
\end{document}